\documentclass[12pt]{amsart}

\usepackage{bbm}
\usepackage{amssymb, amsfonts}
\usepackage[all, arc]{xy}
\usepackage{enumerate}
\usepackage{mathrsfs}
\usepackage{dsfont}
\usepackage{enumitem}
\usepackage{fancyhdr}
\newtheorem{step}{STEP}

\newcommand{\Bp}{B_p}
\newcommand{\Bq}{B_q}
\newcommand{\Sp}{S_p}
\newcommand{\Sq}{S_q}

\newtheorem{theorem}{Theorem}[section]

\theoremstyle{definition}
\newtheorem{definition}[theorem]{Definition}

\newtheorem{notation}[theorem]{Notation}

\theoremstyle{remark}
\newtheorem{remark}[theorem]{Remark}

\newcommand*\diff{\mathop{}\!\mathrm{d}}

\makeatletter
\let\c@equation\c@theorem
\makeatother
\numberwithin{equation}{section}

\title{ Removability Criteria for Separately Harmonic and Subharmonic Functions}

\usepackage[a4paper, margin=1.1in]{geometry}

\author{MANSOUR KALANTAR}

\date{}

\begin{document}
		\maketitle
	\sloppy
	
	\begin{abstract}
Cegrell and Sadullaev showed that a separately harmonic and subharmonic function is subharmonic outside the product of two closed nowhere dense sets. We obtain sufficient conditions, expressed in terms of the second set, under which the function is subharmonic on the whole domain. In particular, we prove that this holds whenever the second set is polar.

	\end{abstract}
	
	\email{Mansour Kalantar// College of San Mateo//kalantarm@smccd.edu}
	
 

	\section{Introduction}

The classical Hartogs theorem \cite{Chabat} states that a function defined on a product of two domains in complex Euclidean spaces and holomorphic in each variable separately is jointly holomorphic. This naturally raises the question whether functions that belong separately to a given class must also belong jointly to that class. For subharmonic functions, the answer is negative: Wiegerinck \cite{Wiegerinck} constructed a separately subharmonic function that is not jointly subharmonic (see also Wiegerinck and Zeinstra \cite{Wiegerinck Zeinstra}).

 An intriguing question, raised by Wiegerinck and Zeinstra \cite{Wiegerinck Zeinstra}, is whether a function \(u(x,y)\), harmonic in \(x\) and subharmonic in \(y\), is necessarily subharmonic. Several authors have obtained affirmative answers under additional assumptions. Arsove \cite{Arsove} assumed that \(u\) is bounded below by a locally integrable function. Imomkulov \cite{Imm} (as cited in \cite{Cegrell Sadullaev}) required that the partial functions \(u(x,\cdot)\) be real analytic, while Kołodziej and Thorbiörnson \cite{Ko-Th} assumed that they belong to \(C^{2}\). These results were subsequently generalized by Riihentaus in a series of papers \cite{Rii2,Rii3} (see also \cite{Riibook} and the references therein). Sadullaev \cite{Sad} continued the study of these functions and established a Riesz decomposition theorem for them. On the other hand, Cegrell and Sadullaev \cite{Cegrell Sadullaev} proved that there exist closed nowhere dense sets \(S_p\) and \(S_q\) in the domains on the product of which \(u\) is defined, such that \(u\) is subharmonic on the complement of \(S_p\times S_q\).

In the present paper, we do not impose any additional assumptions on \(u\) or on its partial functions. Instead, we derive sufficient conditions in terms of one projection of the singular set introduced by Cegrell and Sadullaev. We first show that \(u\) is subharmonic if, for every $x$,  the Riesz measure \(\mu_x\) of the partial function \(u(x,\cdot)\) does not charge \(S_q\). As a consequence, we prove that \(u\) is subharmonic if \(S_q\) is polar.

	\section{Definitions and Preliminaries}
	
	\begin{notation}
	
		By natural numbers, we mean the set
		$$
		\mathbb{N}=\{1,2,3,\dots\}.
		$$
		
		For $\xi \in \mathbb{R}^N$, with ($N \geq 2$), we denote by $B_N(\xi,r)$ the open ball in $\mathbb{R}^N$ centered at $\xi$ and of radius $r>0$. The symbol $\lambda_N$ denotes the $N-$dimensional Lebesgue measure. The closure of a set $E$ is denoted by $\overline{E}$, and its characteristic function is denoted by $\mathbf{1}_E$, that is,
		$$
		\mathbf{1}_E(x)=
		\begin{cases}
			1, & \text{ if }x\in E,\\
			0, & \text{ if } x\notin E.
		\end{cases}
		$$
		
		Given an open set $\Omega \subset \mathbb{R}^N$, we denote by $C_c^\infty(\Omega)$ the space of smooth functions whose support is compact and contained in $\Omega$. The support of a function $\varphi$, denoted by $\operatorname{supp}(\varphi)$, is the closure of the set
		$$\{x\in \Omega : \varphi(x)\neq 0\}.$$
		
		If $\Omega$ is an open subset of $\mathbb{R}^N$, $T$ is a distribution on $\Omega$, and $\varphi \in C_c^\infty(\Omega)$, we write $\langle T,\varphi\rangle$ for the action of $T$ on $\varphi$. The space of distributions on $\Omega$, that is, the space of continuous linear functionals on $C_c^\infty(\Omega)$, is denoted by $\mathcal{D}'(\Omega)$.
		
		We denote by $\Delta_\xi$ the Laplacian with respect to the $\xi-$variable:
		\[
		\Delta_\xi:=\sum_{i=1}^{N}\frac{\partial^2}{\partial \xi_i^2}.
		\]

	\end{notation}
	
	Throughout this paper, we fix two Euclidean balls  \(\Bp\subset \mathbb{R}^p\), and $B_q\subset\mathbb{R}^q$ with $p,q\ge 2$.

	\begin{definition}
		A function \(u(x, y)\) defined on \(B_p \times B_q\) is said to be \emph{separately harmonic and subharmonic} if for each fixed \((a, b) \in B_p \times B_q\), the function \(x \mapsto u(x, b)\) is harmonic on \(B_p\), and \(y \mapsto u(a, y)\) is subharmonic on \(B_q\). If \(u\) is subharmonic in each variable, it is called \emph{separately subharmonic}. If \(u\) is harmonic in each variable, it is called \emph{separately harmonic}.
	\end{definition}
	
	We recall the following results:
	
	\begin{theorem}[Cegrell and Sadullaev \cite{Cegrell Sadullaev}]\label{cegrell sadullaev}
		Let \(u(x, y)\) be separately harmonic-subharmonic on a neighborhood of \(\overline{B_p \times B_q}\). Then there exist two nowhere dense closed sets \(\Sp \subset B_p\) and \(S_q \subset B_q\) such that \(u\) is subharmonic on \(\left((B_p \setminus S_p) \times B_q \right)\cup \left(B_p \times (B_q \setminus S_q)\right)\). The sets $S_p$ and $S_q$ will be called the projections of the singular set of $u$.
	\end{theorem}
	
	\begin{theorem}[Avanissian \cite{Avanissian}]\label{Avaniss}
		Separately subharmonic functions that are locally bounded above are subharmonic.
	\end{theorem}
	
	\begin{theorem}[Lelong \cite{Lelong 2}]\label{Lelong}
		Separately harmonic functions are harmonic.
	\end{theorem}
	
	\begin{theorem}[Generalized Weyl's Lemma]\label{Weyl}
		Let \(\Omega \subset \mathbb{R}^N\) be open and let \(T \in \mathcal{D}'(\Omega)\) be a distribution satisfying
		\[ \langle T, \Delta \phi \rangle = 0 \quad \text{for all } \phi \in C_c^\infty(\Omega). \]
		Then \(T\) is a regular distribution given by a harmonic function.
	\end{theorem}
	See Donoghue \cite[p.~127]{distributions}.

	\medskip
	\noindent
	\textbf{Riesz Decomposition Theorem.}
	Let \(u\) be a subharmonic function in a ball \(B \subset \mathbb{R}^N\), with \(u \not\equiv -\infty\). Then there exists a unique Radon measure \(\mu\) on \(B\), called the \emph{Riesz measure} of \(u\), and a harmonic function \(h\) on \(B\), such that
	\[ u(x) = h(x) - \int_B G_B(x, y) \, d\mu(y), \]
	where \(G_B(x, y)\) is the Green's function of \(B\) with pole at $y=x$. Moreover,
	\[\mu= \Delta u \quad \text{in the distributional sense}. \] (See Armitage and Gardiner \cite[pg 105]{Armitage Gardiner book})
	
	\medskip
	\noindent
	If \(\mu\) is a Radon measure and \(u \in L^1_{\text{loc}}(\Omega)\), we say \(\Delta u = \mu\) in the distributional sense if
	\[ \int_\Omega u \, \Delta \varphi \, d\lambda_N = \int_\Omega \varphi \, d\mu \quad \text{for all } \varphi \in C_c^\infty(\Omega). \]

	\medskip
	\noindent

	\section{Main Results}
\begin{theorem}\label{ThWith_mu=0}
	Let \(u(x,y)\) be a separately harmonic and subharmonic function defined on a neighborhood of $
	\overline{B_p\times B_q},$ and let $\mu_x$ designate the Riesz measure of $u(x,\cdot)$.
	Let \(S_q\subset B_q\) denote the second projection of the singular set of \(u\) given by Theorem 2.3.  If 
	\begin{equation}\label{measureSq0}
		\mu_x(S_q)=0
	\end{equation}
	for all $x\in B_p$, then  \(u\) is subharmonic on \(B_p\times B_q\).
\end{theorem}
	\begin{proof}

Let $(\varepsilon_n)$ be a sequence of positive numbers stricly decreasing to zero. For each natural number $n$, define 
$$U_n:=\{\xi\in \mathbb{R}^q:\operatorname{dist}(\xi,\Sq)<\varepsilon_n\}.$$ 
Then, choose smooth functions $\chi_n \in C_c^\infty(U_n)$ such that $ 0 \leq \chi_n \leq 1,$ and  $\chi_n \equiv 1$ on $\overline{U_{n+1}}$.
The sequence $\chi_n(z)$ decreases pointwise to the characteristic function $\mathbf{1}_{S_q}(z)$ of $\Sq$, as $n\to+\infty$.

Next, without loss of generality, we may assume that $B_q=B_q(0,R)$ is the ball centered at the origin with radius $R>0$. Choose smooth exhaustion functions
$\theta_m \in C_c^\infty(\Bq\!\left(0,R-\varepsilon_{m+1}\right))$ such that	$0 \leq \theta_m \leq 1$ and $\theta_m \equiv 1 $ on $\overline{	\Bq\!\left(0,R-\varepsilon_{m}\right)}.$
The sequence $	(\theta_m)$ increases to the characteristic function $\mathbf{1}_{B_q}$, as $m\to+\infty$.


\begin{step}\label{step rmn harmonic}
	Fix \(y\in B_q\). If necessary, replace $(\varepsilon_m)$ by a tail of the sequence and renumber it so that $ B_q(y,\varepsilon_1)$ is relatively compact in $\Bq$.  If $y\in \Bq\setminus\Sq$, choose smooth functions
	\[
	\beta_m\in C_c^\infty(\Bq(y,\varepsilon_m))
	\]
	such that \(0\le \beta_m\le 1\), and \(\beta_m\equiv 1\) on \(\overline{\Bq(y,\varepsilon_{m+1})}\); if $y\in \Sq$, then set $\beta_m\equiv0$.  For all fixed natural numbers \(n,m\), define
	\[
	r_{n,m}(x,y)
	:=
	\int_{B_q}
	\theta_m(z)
	(1-\chi_n(z))
	(1-\beta_m(z))
	G_{B_q}(y,z)
	\,d\mu_x(z),
	\]where $G_{B_q}$ designates the Green function of the ball $\Bq$.
	Then, the function \(r_{n,m}(\cdot,y)\) is harmonic on \(B_p\).
\end{step}

Set
\[
\eta_{n,m}(y,z)
:=
\theta_m(z)
(1-\chi_n(z))
(1-\beta_m(z))
G_{B_q}(y,z).
\]
If $y\in \Bq\setminus\Sq$, the factor \(1-\beta_m\) vanishes in a neighborhood of \(y\), and therefore removes the singularity of \(G_{B_q}(y,\cdot)\). If $y\in \Sq$, the factor $(1-\chi_n)$ that vanishes on a neighboorhood of $\Sq$, remove the singular point $y$. Moreover, since  the smooth functions $(1-\chi_n)$,  $(1-\beta_m)$ and $G_{\Bq}(y,\cdot)$ are multipled by the test function \(\theta_m\in C_c^\infty(B_q)\), the integrand $\eta_{n,m}$ becomes a test function compactely supported in $B_q$.

To prove harmonicity of $r_{n,m}$, take $\varphi\in C^\infty_c(B_p)$. Then,
\begin{align*}
	\langle \Delta_x r_{n,m}(x,y),\varphi(x)\rangle&=\int_{B_p}r_{n,m}(x,y)\Delta_x\varphi(x)\diff\lambda_p(x)\\
	&=\int_{B_p}\left(\int_{B_q}u(x,z)\Delta_z\eta_{n,m}(y,z)
	\Delta_x\varphi(x)\diff\lambda_q(z)\right)\diff\lambda_p(x).
\end{align*}
where the  second equality follows from the definition of the Riesz measure,
\[
r_{n,m}(x,y)
=
\int_{B_q}
u(x,z)\Delta_z\eta_{n,m}(y,z)\,\diff\lambda_q(z).
\]
The support of \(\Delta_z\eta_{n,m}\), which is included into the support of $\eta_{n,m}$,  is compactly contained in \(B_q\setminus S_q\). Indeed, near \(S_q\), the function \(1-\chi_n\) is identically zero, and near \(y\in \Bq\setminus\Sq\), the function \(1-\beta_m\) is identically zero. Near the boundary of \(B_q\), the cutoff \(\theta_m\) guarantees compact support.  It follows that we may apply Fubini's theorem to the above double integral, and this one equals
$$\int_{B_q}\left(\int_{B_p}u(x,z) \Delta_x\varphi(x)\diff\lambda_p(x)\right)\Delta_z\eta_{n,m}(y,z)
\diff\lambda_q(z).$$
The inner integral equals the laplacian of the harmonic function $x\to u(x,z)$ with respect to $x$ in the sense of distributions, for fixed $z$. It vanishes, since $u(.,z)$ is harmonic by hypothesis, and we obtain 
$$\langle\Delta_x r_{n,m}(x,y),\varphi(x)\rangle=0.$$ Hence $r_{n,m}(.,y)$ is harmonic on $\Bp$.


\begin{step}\label{step r harmonic}
	For fixed \(y \in B_q \), set
	\[
	r(x,y)
	:=
	\int_{B_q\setminus S_q}
	G_{B_q}(y,z)\,d\mu_x(z).
	\]
	Then, the function \(r(\cdot,y)\) is harmonic on \(B_p\).
\end{step}

Let us first prove that \(G_{B_q}(y,\cdot)\) is \(\mu_x\)-integrable, and that the Riesz measure of \(u(x,\cdot)\) does not charge \(\{y\}\). With 
$$p(x,y)=\int_{B_q}G_{B_q}(y,z)\diff\mu_x(z),$$ write
\[
u(x,\cdot)=h(x,\cdot)-p(x,\cdot).
\]
If \(p(x,y)=+\infty\), then
\[
u(x,y)=-\infty.
\]
But this is impossible, because for fixed \(y\), the function \(u(\cdot,y)\) is harmonic on \(B_p\), and hence finite. Therefore
\begin{equation}\label{p is integrable}
	p(x,y)<+\infty.
\end{equation}
In particular, \(G_{B_q}(y,\cdot)\) is integrable with respect to the measure \(\mu_x\). 

Moreover, we have
$$0\le \int_{\{y\}}G_{\Bq}(y,z)\,d\mu_x(z) \le p(x,y) .  $$
If the measure $\mu_x$ charged $\{y\}$ the above integral would be $+\infty$, since \(G_{B_q}(y,\cdot)\) has a pole at \(y\).  The finiteness of \(p(x,y)\) implies
\begin{equation}\label{p is finite}
	\mu_x(\{y\})=0.
\end{equation}

Next,  let \(m\to\infty\). Since
\[
0
\leq
\theta_m(z)
(1-\chi_n(z))
(1-\beta_m(z))
G_{B_q}(y,z)
\leq
G_{B_q}(y,z),
\]
and both \(\theta_m\) and \(1-\beta_m\) increase pointwise to $1$ and $\mathbf{1}_{B_q\setminus\{y\}}$
respectively, the monotone convergence theorem\footnote{Or the dominated convergence theorem and $\mu_x-$integrability of $G_{\Bq}(y,.)$, due to (\ref{p is integrable}).} yields
\[
r_{n,m}(x,y)\to r_n(x,y),
\]
where
\[
r_n(x,y)
:=
\int_{B_q}
(1-\chi_n(z))
G_{B_q}(y,z)\,d\mu_x(z).
\]
In view of (\ref{p is finite}), removing the point \(y\) from \(B_q\) does not affect the integral. Hence, by Harnack's theorem, the function \(r_n(\cdot,y)\) is either identically \(+\infty\), or harmonic on \(B_p\). Since
\[
0\le r_n(x,y)\le p(x,y)<+\infty,
\]
the first possibility is excluded, and so \(r_n(\cdot,y)\) is harmonic on \(B_p\).

Finally,  let \(n\to\infty\). Since \(1-\chi_n\) increases pointwise to \(\mathbf{1}_{B_q\setminus S_q}\), another application of the monotone convergence theorem yields
\[
r_n(x,y)\to r(x,y).
\]
Since
\[
0\le r(x,y)\le p(x,y)<+\infty,
\]
Harnack's theorem implies that the function  \(r(\cdot,y)\) is not identically $+\infty$, and so is harmonic on \(B_p\).


\begin{step}
	For each fixed $x\in B_p$, let
	\[
	u(x,y)
	=
	h(x,y)-p(x,y)
	\]
	be the Riesz decomposition of the subharmonic function
	$u(x,\cdot)$ on $B_q$, where $h(x,\cdot)$ is harmonic on $B_q$, and
	\[
	p(x,y)
	=
	\int_{B_q}
	G_{B_q}(y,z)\,d\mu_x(z).
	\]
	Then the function $h(x,y)$ is harmonic on $B_p\times B_q$ and $u$ is
	subharmonic on $B_p\times B_q$.
\end{step}

Notice that for every fixed $y\in B_q$, the potential in the Riesz decomposition can be written as
\[
p(\cdot,y)
=
q(\cdot,y)
+
r(\cdot,y),
\]
where $$q(x,y)=\int_{S_q} G_{B_q}(y,z)d\mu_x(z).$$  
By (\ref{measureSq0}),  we have $q(x,y)=0$ and so $p(x,y)=r(x,y).$
By Step \ref{step r harmonic}, the function \(r(\cdot,y)\) is  harmonic on \(B_p\). Thus, so is \(p(\cdot,y)\). Since \(u(\cdot,y)\) is harmonic by hypothesis, it follows that
\[
h(\cdot,y)
=
u(\cdot,y)
+
p(\cdot,y)
\]
is harmonic on \(B_p\) for every fixed $y\in B_q.$

On the other hand, for every fixed \(x\in B_p\), the second partial function $h(x,\cdot)$
is harmonic on \(B_q\) by construction. Thus $h(x,y)$ is separately harmonic on $\Bp\times\Bq$, and by Lolong's theorem, harmonic on $\Bp\times\Bq$.

Finally,
\[
u=h-p\le h,
\]
since \(p\) is nonnegative. As \(h\) is harmonic on $B_p\times B_q,$ it is locally bounded above. Hence \(u\) is also locally bounded above on $B_p\times B_q.$
Since \(u\) is separately subharmonic, Avanissian's theorem applies, and we conclude that \(u\) is subharmonic on $B_p\times B_q.$
	
	\end{proof}
	
\begin{remark}
	The conclusion of Theorem \ref{ThWith_mu=0} remains valid if condition (\ref{measureSq0}) is replaced by the following slightly weaker assumptions:
	\begin{enumerate}\label{q>0}
		\item[(i)] $\mu_x(S_q)=0$ for all $x\in S_p$;
		\item[(ii)] for each connected component $D$ of $\Bp\setminus\Sp$, there exists a point $x_D\in D$ such that $\mu_{x_D}(S_q)=0$.
	\end{enumerate}
\end{remark}

\begin{proof}
For every nonnegative $\psi\in C_c^\infty(B_q)$ define
	\[
	q_{n,\psi}(x):=\int_{B_q}\chi_n(z)\psi(z)\,d\mu_x(z)
	=\int_{B_q}u(x,z)\Delta_z\bigl(\chi_n(z)\psi(z)\bigr)\,dz.
	\]
	Exactly as in Step 1 of the proof of Theorem \ref{ThWith_mu=0}, using that for every $\varphi\in C_c^\infty(B_p\setminus S_p)$ the support of
	$\Delta_x\varphi(x)\,\Delta_z(\chi_n\psi)$ is contained in $(B_p\setminus S_p)\times B_q$, where $u$ is locally integrable by Theorem \ref{cegrell sadullaev}, one proves that $q_{n,\psi}$ is harmonic on $B_p\setminus S_p$.
	
	Since $\chi_n\downarrow\mathbf 1_{S_q}$, the dominated convergence theorem implies that $q_{n,\psi}(x)\to q_{\psi}(x)$, where
	\[
	q_{\psi}(x):=\int_{S_q}\psi(z)\,d\mu_x(z).
	\]
	As each $q_{n,\psi}$ is nonnegative and harmonic, Harnack's convergence theorem implies that $q_{\psi}$ is harmonic on $B_p\setminus S_p$.
	
	Since $\psi\ge0$, we have $q_{\psi}\ge0$ on $B_p\setminus S_p$. Hence, by the strong maximum principle, on each connected component $D$ of $B_p\setminus S_p$ either $q_{\psi}\equiv0$ or $q_{\psi}>0$. The latter alternative is excluded by assumption (ii). Since this argument holds for every nonnegative $\psi\in C_c^\infty(B_q)$, we conclude that $\mu_x(S_q)=0$ for all $x\in B_p\setminus S_p$. Together with assumption (i), this yields condition (\ref{measureSq0}), and the conclusion follows from Theorem \ref{ThWith_mu=0}.
\end{proof}

\begin{theorem}
	Let \(u(x,y)\) be a separately harmonic and subharmonic function defined on a neighborhood of $
\overline{B_p}\times \overline{B_q}.$
Let \(S_q\subset B_q\) denote the second projection of the singular set of \(u\) given by Theorem 2.3. Suppose that \(S_q\) is polar. Then \(u\) is subharmonic on \(B_p\times B_q\).
\end{theorem}

\begin{proof}
Since \(S_q\) is polar, it has Newtonian capacity zero. In view of the previous proposition, it is sufficient to show that
\begin{equation}\label{mu_x(S_q)=0}
	\mu_x(S_q)=0.
\end{equation}

Suppose, toward a contradiction, that \(\mu_x(S_q)>0\). Let
\[
p(x,y):=\int_{B_q}G_{B_q}(y,z)\,d\mu_x(z).
\]By  (\ref{p is integrable}) we have  \(p(x,y)<+\infty\) for every \(y\in B_q\). Let
$$E_m:=S_q\cap \{y\in B_q:p(x,y)\le m\}.$$  Then
\[
S_q=\bigcup_{m=1}^{\infty}
E_m.
\]
Therefore, for some set $E_m$, we must have  \(\mu_x(E_m)>0\).  
By inner regularity, we may choose a compact set \(K \subset E_m\) such that
\[
\mu_x(K)>0,
\]
and let $
\nu:=\mu_x|_K$
be the restriction of \(\mu_x\) to \(K\).

We estimate the Green energy of \(\nu\). Since \(G_{B_q}\geq 0\),
\[
\begin{aligned}
	I_{B_q}(\nu)
	&:=
	\int_{B_q} G_{B_q}\nu(y)\, d\nu(y) \\
	&=
	\int_K
	\left(
	\int_K G_{B_q}(y,z)\, d\nu(z)
	\right)
	d\nu(y) \\
	&\leq
	\int_K
	\left(
	\int_{B_q} G_{B_q}(y,z)\, d\mu_x(z)
	\right)
	d\nu(y).
\end{aligned}
\]
By definition of the potential \(p(x,\cdot)\),
\[
\int_{B_q} G_{B_q}(y,z)\, d\mu_x(z)
=
p(x,y).
\]
Since \(K \subset E_m\), we have \(p(x,y)\leq m\) for every
\(y\in K\), and therefore
\[
I_{B_q}(\nu)
\leq
\int_K p(x,y)\, d\nu(y)
\leq
m\,\nu(K)
<
+\infty.
\]

We have thus obtained a nonzero positive measure $\nu$, supported on the compact polar set $K$, and having finite energy. This contradicts the characterization of compact polar sets: a compact set is polar if and only if it is null for every positive measure of finite energy and supported on it (for a slightly more general result, see Helms \cite[Theorem 6.2.3]{Helms}). This contradiction proves \eqref{mu_x(S_q)=0}.

\end{proof}

	\section*{Declarations}
	
	\subsection*{Ethical Approval}
	Not applicable.
	
	\subsection*{Funding}
	The author received no funding for this work.
	
	\subsection*{Availability of data and materials}
	Not applicable. This manuscript does not report data generation or analysis.


\begin{thebibliography}{99}
		
		\bibitem{Armitage Gardiner book} Armitage, D.M. and Gardiner, S.J. \emph{Classical potential theory}, Springer-Verlag, London, 2001.
	
		\bibitem{Arsove} Arsove, M.G. \emph{On the subharmonicity of double subharmonic functions}, Proc. Amer. Math. Soc. 17, pp. 622-626, 1966.
		
		\bibitem{Avanissian} Avanissia, V. \emph{Fonctions plurisousharmoniques et fonctions doublement sousharmoniques},
		Ann. Sci. Ecol. Norm. Super. 78, pp. 101-161, 1961.
		
		\bibitem{Cegrell Sadullaev} Cegrell, U. and Sadullaev, A. \emph{Separately subharmonic functions }, Uzbek. Math. J. 1, pp. 78-83, 1993.
	
	
		\bibitem{Chabat} Chabat, B.  \emph{ Introduction à l’analyse complexe}, Tome 2: Fonctions de plusieurs variables, Éditions Mir, Moscow, 1990.
	
		\bibitem{distributions} Donoghue Jr, W.F. \emph{Distributions and Fourier Transforms }, Pure and Applied Mathematics, Volume 32, 1969.
		
		\bibitem{Helms} Helms, L. L.
	\emph{Potential Theory},
	Universitext, Springer, London, 2009.
	
		\bibitem{Imm}Imomkulov, S.A. \emph{Separately Subharmonic functions}, (in Russian), Dokl. USSR. 2l, pp. 8-10, 1990.
		
		
		
		\bibitem{Ko-Th} Kołodziej, S. and Thorbi\"{o}rnson, J. \emph{Separately harmonic and subharmonic functions}, Potential Anal. 5, pp. 463-466, 1996.
		
		\bibitem{Lelong 2}Lelong, P. \emph{Fonctions plurisousharmoniques et fonctions analytiques de variables r\''{e}elles}, Ann. Inst. Fourier 11, pp. 515-562, 1961.
	
		\bibitem{Rii2} Riihentaus, J. \emph{On separately subharmonic and harmonic functions}, Int. J. Pure Appl. Math., 35, No. 4, pp. 435 -446, 2007.
		
		\bibitem{Rii3} Riihentaus, J. \emph{Subharmonic functions, generalizations and separately subharmonic functions}, in: The XIV-th Conference on Analytic Functions, Chelm, Section of Mathematics and Computer Science, 2, pp 49-76, 2007.
		\bibitem{Riibook} Riihentau, J. \emph{Subharmonic functions, generalizations, holomorphic Functions, meromorphic Functions, and properties}, Bentham Books, 2021.
		
		\bibitem{Sad} Sadullaev, A. \emph{On separately subharmonic functions (Lelong''s problem)}, Annales de la Faculté des Sciences de Toulouse, Vol. XX, n$^\circ$ Spécial, pp. 183-187, 2011.
		
		
		\bibitem{Wiegerinck} Wiegerinck, J. \emph{Separately subharmonic functions need not be subharmonic}, Proc. Amer. Math. Soc. 104, pp. 770-771, 1988.
		
		\bibitem{Wiegerinck Zeinstra} Wiegerinck, J. and Zeinstra, R. \emph{Separately subharmonic functions: when are they subharmonic?}, Proc. of Symposia in Pure Math. 52, pp. 245-249, 1991 part 1.
	\end{thebibliography}
\end{document}